\documentclass{amsproc}
\usepackage{amssymb}
\usepackage{amsfonts}

\theoremstyle{plain}

\newtheorem{remark}{Remark}

\newtheorem{theorem}{Theorem}
\input{tcilatex}

\begin{document}
\title{On Boman's Theorem On Partial Regularity Of Mappings}
\author{Tejinder S. Neelon}
\address{Department of Mathematics\\
California State University San Marcos\\
San Marcos, CA 92096-0001, USA}
\email{neelon@csusm.edu}
\urladdr{http://www.csusm.edu/neelon/neelon.html}

\begin{abstract}
Let $\Lambda \subset \mathbb{R}^{n}\times \mathbb{R}^{m}$ and $k$ be a
positive integer. Let $f:\mathbb{R}^{n}\rightarrow \mathbb{R}^{m}$ be a
locally bounded map such that for each $(\xi ,\eta )\in \Lambda ,$ the
derivatives $D_{\xi }^{j}f(x):=\left. \frac{d^{j}}{dt^{j}}f(x+t\xi
)\right\vert _{t=0},$ $j=1,2,...k,$\ exist and are continuous$.$ In order to
conclude that any such map $f$ is necessarily of class $C^{k}$ it is
necessary and sufficient that $\Lambda $ be \emph{not} contained in the
zero-set of a nonzero homogenous polynomial $\Phi (\xi ,\eta )$ which is
linear in $\eta =(\eta _{1},\eta _{2},...,\eta _{m})$ and homogeneous of
degree $k$ in $\xi =(\xi _{1},\xi _{2},...,\xi _{n})$.

This generalizes a result of J. Boman for the case $k=1$. The statement and
the proof of a theorem of Boman for the case $k=\infty $ is also extended to
include the Carleman classes $C\{M_{k}\}$ and the Beurling classes $%
C(M_{k}). $(\cite{Boman-partial})
\end{abstract}

\maketitle

A continuous function $f:%
%TCIMACRO{\U{211d} }%
%BeginExpansion
\mathbb{R}
%EndExpansion
^{n}\rightarrow 
%TCIMACRO{\U{211d} }%
%BeginExpansion
\mathbb{R}
%EndExpansion
$\ that is differentiable when restricted to arbitrary differentiable curves
is not necessarily differentiable as a function of several variables \cite%
{rudin}.\textrm{\ }Indeed, there are discontinuous functions $f:%
%TCIMACRO{\U{211d} }%
%BeginExpansion
\mathbb{R}
%EndExpansion
^{n}\rightarrow 
%TCIMACRO{\U{211d} }%
%BeginExpansion
\mathbb{R}
%EndExpansion
$\ whose restrictions to arbitrary analytic arcs are analytic \cite{BMP}.%
\textrm{\ }But a $C^{\infty }$\ function $f:%
%TCIMACRO{\U{211d} }%
%BeginExpansion
\mathbb{R}
%EndExpansion
^{n}\rightarrow 
%TCIMACRO{\U{211d} }%
%BeginExpansion
\mathbb{R}
%EndExpansion
$\ whose restriction to every line segment is real analytic is necessarily
real analytic (\cite{Siciak}). In \cite{Neelon1}, \cite{N2} , \cite{n3} and 
\cite{n4} this result was extended by considering restrictions to algebraic
curves and surfaces of functions belonging to more general classes of
infinitely differentiable functions. It is also well known that a function $%
f:%
%TCIMACRO{\U{211d} }%
%BeginExpansion
\mathbb{R}
%EndExpansion
^{n}\rightarrow 
%TCIMACRO{\U{211d} }%
%BeginExpansion
\mathbb{R}
%EndExpansion
$\ that is infinitely differentiable in each variable separately may be no
better than measurable (\cite{krantz}).\textrm{\ }In \cite{Boman-partial},
the obverse problem is considered; for vector valued functions hypothesis is
made on the source as well as the target space. In this note, Theorem 4 of 
\cite{Boman-partial} is generalized to $C^{k},$ $k\geq 1,$ the class of
functions that have continuous derivatives up to order $k$.

Let $f:%
%TCIMACRO{\U{211d} }%
%BeginExpansion
\mathbb{R}
%EndExpansion
^{n}\rightarrow 
%TCIMACRO{\U{211d} }%
%BeginExpansion
\mathbb{R}
%EndExpansion
^{m}$ be a locally bounded map. For $(\xi ,\eta )\in 
%TCIMACRO{\U{211d} }%
%BeginExpansion
\mathbb{R}
%EndExpansion
^{n}\times 
%TCIMACRO{\U{211d} }%
%BeginExpansion
\mathbb{R}
%EndExpansion
^{m},$\ set\ 
\begin{equation*}
D_{\xi }\left\langle f,\eta \right\rangle (x):=\left. \frac{d}{dt}%
\left\langle f(x+t\xi ),\eta \right\rangle \right\vert _{t=0}\text{ in the
sense of distributions,}
\end{equation*}%
where $\left\langle \cdot ,\cdot \right\rangle $\ denotes the inner product
on $%
%TCIMACRO{\U{211d} }%
%BeginExpansion
\mathbb{R}
%EndExpansion
^{m}.$\ By the Leibniz Integral rule, we have\ 
\begin{equation*}
\frac{d}{dt}\int \left\langle f(x+t\xi ),\eta \right\rangle dx=\int \frac{d}{%
dt}\left\langle f(x+t\xi ),\eta \right\rangle dx\text{.}
\end{equation*}

Let $k,1\leq k<\infty ,$\ be fixed.\ For $\xi \in 
%TCIMACRO{\U{211d} }%
%BeginExpansion
\mathbb{R}
%EndExpansion
^{n},$\ denote by $C_{\xi }^{k}\left( 
%TCIMACRO{\U{211d} }%
%BeginExpansion
\mathbb{R}
%EndExpansion
^{n}\right) $\ the space of all continuous functions $f:%
%TCIMACRO{\U{211d} }%
%BeginExpansion
\mathbb{R}
%EndExpansion
^{n}\rightarrow 
%TCIMACRO{\U{211d} }%
%BeginExpansion
\mathbb{R}
%EndExpansion
$\ such that the derivatives $D_{\xi }^{j}f(x):=\left. \frac{d^{j}}{dt^{j}}%
f(x+t\xi )\right\vert _{t=0},$ $j=1,2,...k,$\ exist and are continuous.\
Similarly, $C_{\xi }^{\infty }\left( 
%TCIMACRO{\U{211d} }%
%BeginExpansion
\mathbb{R}
%EndExpansion
^{n}\right) :=\cap _{k=0}^{\infty }C_{\xi }^{k}\left( 
%TCIMACRO{\U{211d} }%
%BeginExpansion
\mathbb{R}
%EndExpansion
^{n}\right) .$

We are interested in finding the necessary and sufficient conditions on a
subset $\Lambda \subset 
%TCIMACRO{\U{211d} }%
%BeginExpansion
\mathbb{R}
%EndExpansion
^{n}\times 
%TCIMACRO{\U{211d} }%
%BeginExpansion
\mathbb{R}
%EndExpansion
^{m}$\ to have the following property: 
\begin{eqnarray*}
\text{if }f &:&%
%TCIMACRO{\U{211d} }%
%BeginExpansion
\mathbb{R}
%EndExpansion
^{n}\rightarrow 
%TCIMACRO{\U{211d} }%
%BeginExpansion
\mathbb{R}
%EndExpansion
^{m}\text{ is locally bounded } \\
\text{ such that }\left\langle f,\eta \right\rangle &\in &C_{\xi }^{k}\left( 
%TCIMACRO{\U{211d} }%
%BeginExpansion
\mathbb{R}
%EndExpansion
^{n}\right) ,\forall (\xi ,\eta )\in \Lambda ,\text{ then }f\in C^{k}\left( 
%TCIMACRO{\U{211d} }%
%BeginExpansion
\mathbb{R}
%EndExpansion
^{n}\right) .
\end{eqnarray*}

The case $k=1$\ and $k=\infty $ was dealt in \cite{Boman-partial}. \ 

Let $%
%TCIMACRO{\U{2124} }%
%BeginExpansion
\mathbb{Z}
%EndExpansion
_{+}^{n}$\ denote all $n$-tuples of nonnegative integers.\textrm{\ } For $%
\alpha =\left( \alpha _{1},\alpha _{2},...,\alpha _{n}\right) \in 
%TCIMACRO{\U{2124} }%
%BeginExpansion
\mathbb{Z}
%EndExpansion
_{+}^{n},$\ set $|\alpha |=\alpha _{1}+\alpha _{2}+...+\alpha _{n}.$The set $%
%TCIMACRO{\U{2124} }%
%BeginExpansion
\mathbb{Z}
%EndExpansion
_{+}^{n}$ of multi-indices \ is assumed to be ordered lexicographically i.e.
for\ $\alpha =\left( \alpha _{1},\alpha _{2},...,\alpha _{n}\right) ,\beta
=\left( \beta _{1},\beta _{2},...,\beta _{n}\right) \in 
%TCIMACRO{\U{2124} }%
%BeginExpansion
\mathbb{Z}
%EndExpansion
_{+}^{n},$\ define\ $\alpha \prec \beta $\ if there is $i,1\leq i\leq n,$\
such that $\alpha _{1}=\beta _{1},\alpha _{2}=\beta _{2},...,\alpha
_{i-1}=\beta _{i-1},\alpha _{i}<\beta _{i}$.

Let $k_{n}=\binom{k+n-1}{k}$\ denote the number of monomials of degree $k$\
in $n$\ variables.

Then for any $\varphi \in C_{c}^{\infty }(%
%TCIMACRO{\U{211d} }%
%BeginExpansion
\mathbb{R}
%EndExpansion
^{n}),$\ we have\textrm{\ }%
\begin{equation*}
\int D_{\xi }\left\langle f,\eta \right\rangle (x)\varphi (x)dx=\left. \frac{%
d}{dt}\int \left\langle f(x+t\xi ),\eta \right\rangle \varphi
(x)dx\right\vert _{t=0}
\end{equation*}%
\begin{equation*}
=\left. \frac{d}{dt}\left\langle \int f(x)\varphi (x-t\xi )dx,\eta
\right\rangle \right\vert _{t=0}=\left. \left\langle \int f(x)\frac{d}{dt}%
\varphi (x-t\xi )dx,\eta \right\rangle \right\vert _{t=0}
\end{equation*}%
\begin{equation*}
=-\sum_{i}\xi _{i}\left. \left\langle \int f(x)\partial _{i}\varphi (x-t\xi
)dx,\eta \right\rangle \right\vert _{t=0}=\sum_{i,j}\xi _{i}\eta _{j}\int
\partial _{i}f_{j}(x)\varphi (x)dx.
\end{equation*}

By iteration, we obtain the formula for higher-order distributional
derivatives: 
\begin{equation}
D_{\xi }^{p}\left\langle f,\eta \right\rangle (x)=\sum_{\left\vert \alpha
\right\vert =p}\sum_{j=1}^{m}\xi ^{\alpha }\eta _{j}\partial ^{\alpha
}f_{j}(x).  \label{D^k}
\end{equation}%
Let%
\begin{equation*}
\mathcal{B}_{k}:=\left\{ \Phi (\xi ,\eta )=\sum_{j=1}^{m}\sum_{\left\vert
\alpha \right\vert =k}\varphi _{\alpha j}\xi ^{\alpha }\eta _{j}:\varphi
_{\alpha j}\in 
%TCIMACRO{\U{211d} }%
%BeginExpansion
\mathbb{R}
%EndExpansion
,\alpha \in 
%TCIMACRO{\U{2124} }%
%BeginExpansion
\mathbb{Z}
%EndExpansion
_{+}^{n},j\in 
%TCIMACRO{\U{2124} }%
%BeginExpansion
\mathbb{Z}
%EndExpansion
_{+}\right\} .
\end{equation*}

For any function $\Phi (\xi ,\eta ),$\ set $\left\Vert \Phi \right\Vert
:=\max_{\left\Vert \xi \right\Vert \leq 1,\left\Vert \eta \right\Vert \leq
1}\left\vert \Phi (\xi ,\eta )\right\vert .$\textrm{\ }For a subset $%
K\subset \subset \Lambda ,$\ ($\subset \subset $ denotes the compact
inclusion) put $\left\Vert \Phi \right\Vert _{K}:=\max_{(\xi ,\eta )\in
K}\left\vert \Phi (\xi ,\eta )\right\vert .$

\begin{theorem}
\label{main}Let $\Lambda \subset 
%TCIMACRO{\U{211d} }%
%BeginExpansion
\mathbb{R}
%EndExpansion
^{n}\times 
%TCIMACRO{\U{211d} }%
%BeginExpansion
\mathbb{R}
%EndExpansion
^{m}$\ be a subset and $k$\ be a positive integer.\ The following conditions
are equivalent.

(i) $\Lambda $\ is \emph{not} contained in an algebraic hypersurface defined
by an element of $\mathcal{B}_{k}$\ i.e. 
\begin{equation*}
\Phi \in \mathcal{B}_{k},\left. \Phi \right\vert _{\Lambda }\equiv
0\Rightarrow \Phi \equiv 0;
\end{equation*}%
(ii) there exists a set consisting of $m\cdot k_{n}$\ points 
\begin{equation*}
\left( \xi ^{\ast },\eta ^{\ast }\right) =\left\{ \left( \xi ^{(p)},\eta
^{(p)}\right) \in \Lambda ,p=1,2,...,mk_{n}\right\} \ \text{such that }\det
\Delta \left( \xi ^{\ast },\eta ^{\ast }\right) \neq 0,
\end{equation*}%
\ where 
\begin{equation*}
\Delta \left( \xi ^{\ast },\eta ^{\ast }\right) :=\ \left[ \left( \xi
^{(p)}\right) ^{\alpha }\eta _{j}^{(p)}\right] _{|\alpha |=k,1\leq j\leq
m,1\leq p\leq mk_{n}};
\end{equation*}%
(iii) if $f:%
%TCIMACRO{\U{211d} }%
%BeginExpansion
\mathbb{R}
%EndExpansion
^{n}\rightarrow 
%TCIMACRO{\U{211d} }%
%BeginExpansion
\mathbb{R}
%EndExpansion
^{m}$ is locally bounded and $\left\langle f,\eta \right\rangle \in C_{\xi
}^{k},\forall (\xi ,\eta )\in \Lambda ,$ then $f\in C^{k}\left( 
%TCIMACRO{\U{211d} }%
%BeginExpansion
\mathbb{R}
%EndExpansion
^{n},%
%TCIMACRO{\U{211d} }%
%BeginExpansion
\mathbb{R}
%EndExpansion
^{m}\right) .$

If any one of the above equivalent conditions is satisfied, then there
exists a constant $B$\ depending only on $\Lambda $\ such that the following
inequality holds for all locally bounded maps $f:%
%TCIMACRO{\U{211d} }%
%BeginExpansion
\mathbb{R}
%EndExpansion
^{n}\rightarrow 
%TCIMACRO{\U{211d} }%
%BeginExpansion
\mathbb{R}
%EndExpansion
^{m}:$ 
\begin{equation}
\max_{1\leq j\leq m}\max_{|\alpha |=k}\left\vert \partial ^{\alpha
}f_{j}(x)\right\vert \leq B\cdot \sup_{\left( \xi ,\eta \right) \in \Lambda
}\left\vert D_{\xi }^{k}\left\langle f,\eta \right\rangle (x)\right\vert
,\forall x\in 
%TCIMACRO{\U{211d} }%
%BeginExpansion
\mathbb{R}
%EndExpansion
^{n}.  \label{B}
\end{equation}
\end{theorem}

\begin{proof}
We will prove $(i)\Rightarrow (ii)\Rightarrow (iii)\Rightarrow (i).$

$(i)\Rightarrow (ii).$

Suppose $\det \Delta \left( \xi ^{\ast },\eta ^{\ast }\right) =0$\ for every
set of $mk_{n}$\ elements

$\left( \xi ^{\ast },\eta ^{\ast }\right) =\left\{ \left( \xi ^{(p)},\eta
^{(p)}\right) \right\} _{1\leq p\leq mk_{n}}$ in $\Lambda .$ Fix one such
set $\left( \xi ^{\ast },\eta ^{\ast }\right) $ so that the rank $l:=%
\limfunc{rank}\Delta \left( \xi ^{\ast },\eta ^{\ast }\right) $\ is
positive. \ Let $\Delta ^{(l)}$\ denote some $l\times l$ submatrix of $%
\Delta \left( \xi ^{\ast },\eta ^{\ast }\right) $ such that the minor $\det
\Delta ^{(l)}$ is nonzero. \textrm{\ }Let\textrm{\ }$\Delta ^{(l+1)}$ be a $%
(l+1)\times (l+1)$ submatrix of $\Delta \left( \xi ^{\ast },\eta ^{\ast
}\right) $ that contains $\Delta ^{(l)}$ as a submatrix. Replace the point $%
\left( \xi ^{(p_{0})},\eta ^{(p_{0})}\right) $ in $\Delta ^{(l+1)}$ which
does not appear in $\Delta ^{(l)}$ by variables $(\xi ,\eta )\in 
%TCIMACRO{\U{211d} }%
%BeginExpansion
\mathbb{R}
%EndExpansion
^{n}\times 
%TCIMACRO{\U{211d} }%
%BeginExpansion
\mathbb{R}
%EndExpansion
^{m}$. By expanding $\Delta ^{(l+1)}$ along the row where the replacement
took place we obtain an element\ 
\begin{equation*}
\Phi (\xi ,\eta )=\sum_{\alpha ,j}\varphi _{\alpha j}\xi ^{\alpha }\eta _{j},%
\text{ }
\end{equation*}%
of $\mathcal{B}_{k}$ which is nonzero since one of its coefficients
coincides with $\det $ $\Delta ^{(l)}$ up to a sign.

Since $\Delta \left( \xi ^{\ast },\eta ^{\ast }\right) $ has rank $l,$ we
find that $\Phi (\xi ,\eta )=0$ for all $(\xi ,\eta )\in \left( \xi ^{\ast
},\eta ^{\ast }\right) .$ If $\Phi (\xi ,\eta )=0$ for all $(\xi ,\eta )\in
\Lambda ,$ we are done. Otherwise, choose a point $(\widetilde{\xi },%
\widetilde{\eta })\in \Lambda \smallsetminus \left( \xi ^{\ast },\eta ^{\ast
}\right) $ with $\Phi (\widetilde{\xi },\widetilde{\eta })\neq 0.$

Let $\left( \widetilde{\xi ^{\ast }},\widetilde{\eta ^{\ast }}\right) $ be
the set which is obtained from $\left( \xi ^{\ast },\eta ^{\ast }\right) $
by replacing the point $\left( \xi ^{(p_{0})},\eta ^{(p_{0})}\right) $ by $(%
\widetilde{\xi },\widetilde{\eta }).$ Then, the $\limfunc{rank}\Delta \left( 
\widetilde{\xi ^{\ast }},\widetilde{\eta ^{\ast }}\right) \geq l+1.$ By
repeating above procedure, we find a sequence of subsets $\left( \xi ^{\ast
},\eta ^{\ast }\right) ^{(i)}\subset \Lambda ,$ $i=1,2,3,..,$ each with $%
mk_{n}$\ elements such that the $\limfunc{rank}\Delta \left( \xi ^{\ast
},\eta ^{\ast }\right) ^{(j)}\ $is a strictly increasing sequence of
nonnegative integers. After finitely many steps we obtain a nonzero element
of $\mathcal{B}_{k}$ which vanishes on the entire $\Lambda .$

$(ii)\Rightarrow (iii).$

Let $\left( \xi ^{\ast },\eta ^{\ast }\right) =\left\{ \left( \xi
^{(p)},\eta ^{(p)}\right) \in \Lambda \right\} _{1\leq p\leq mk_{n}}$\ be a
set of points such that

$\det \Delta \left( \xi ^{\ast },\eta ^{\ast }\right) \neq 0.$\textrm{\ }By
applying Cramer's rule to (\ref{D^k}), we get%
\begin{equation*}
\partial ^{\alpha }f_{j}(x)=\sum_{p=1}^{mk_{n}}\frac{\det \Delta _{\alpha
j}^{(p)}}{\det \Delta }D_{\xi ^{(p)}}^{k}\left\langle f,\eta
^{(p)}\right\rangle (x)\text{ in the distributional sense,}
\end{equation*}%
where $\Delta _{\alpha j}^{(p)}$\ denotes the cofactor obtained by deleting
the $(\alpha ,j)$-th row and the $p$-th column.\textrm{\ }Since $D_{\xi
}^{k}\left\langle f,\eta \right\rangle \in C^{0}$ for all $(\xi ,\eta )\in
\Lambda ,$ we have\textrm{\ }%
\begin{equation*}
\partial ^{\alpha }f_{j}(x)=\sum_{p=1}^{mk_{n}}\frac{\det \Delta _{\alpha
j}^{(p)}}{\det \Delta }D_{\xi ^{(p)}}^{k}\left\langle f,\eta
^{(p)}\right\rangle (x)\in C^{0}.\mathrm{\ }
\end{equation*}%
Furthermore, there exists a constant $B=B(k,f,\Lambda )$\ such that 
\begin{equation*}
\left\vert \partial ^{\alpha }f_{j}(x)\right\vert \leq
\sum_{p=1}^{mk_{n}}\left\vert \frac{\det \Delta _{\alpha j}^{(p)}}{\det
\Delta }\right\vert \left\vert D_{\xi ^{(p)}}^{k}\left\langle f,\eta
^{(p)}\right\rangle (x)\right\vert \leq B\cdot \sup_{(\xi ,\eta )\in \Lambda
}\left\vert D_{\xi }^{k}\left\langle f,\eta \right\rangle (x)\right\vert ,
\end{equation*}%
for all $\alpha $ with $|\alpha |=k,$ and all $j=1,2,...,m.$

$(iii)\Rightarrow (i).$

Suppose $(i)$\ does not hold. Let $\Phi \in \mathcal{B}_{k}$\ be such that $%
\left. \Phi \right\vert _{\Lambda }\equiv 0.$\textrm{\ }We can write $\Phi
(\xi ,\eta )=\left\langle \varphi _{\cdot }(\xi ),\eta \right\rangle ,$\
where $\varphi _{\cdot }(\xi ):=\left( \varphi _{1}(\xi ),\varphi _{2}(\xi
),...,\varphi _{m}(\xi )\right) $\ and $\varphi _{j}(\xi )=\sum_{\left\vert
\alpha \right\vert =k}\varphi _{\alpha j}\xi ^{\alpha },j=1,2,...,m,$\
homogeneous polynomials of degree $k.$

Define the map 
\begin{equation*}
f(x):=\left\{ 
\begin{array}{c}
\left( \ln \left\vert \ln \left\vert x\right\vert \right\vert \right)
\varphi _{\cdot }(x)\text{ if }x\neq 0, \\ 
0\text{ if }x=0%
\end{array}%
\right. .
\end{equation*}%
Clearly $f\notin C^{k}$ and $f$ is $C^{\infty }$ in $\{x\in 
%TCIMACRO{\U{211d} }%
%BeginExpansion
\mathbb{R}
%EndExpansion
^{n}:0<|x|<1\}.$ We will prove that $D_{\xi }^{k}\left\langle f(x),\eta
\right\rangle $ exists at $x=0,$ for all $(\xi ,\eta )\in \Lambda .$ It is
easy to see that here are constants $C_{\alpha }$ such that 
\begin{equation*}
\left\vert \partial ^{\alpha }\ln \left\vert \ln \left\vert x\right\vert
\right\vert \right\vert \leq \frac{C_{\alpha }}{|x|^{|\alpha |}\left\vert
\ln |x|\right\vert },\forall \alpha ,|\alpha |\geq 1.
\end{equation*}

Since the $\varphi _{j}(x)$'s are homogeneous polynomials of degree $k,$
when the Leibniz's formula is applied to the products $\left( \ln \left\vert
\ln \left\vert x\right\vert \right\vert \right) \varphi _{j}(x)$, it is
clear that all terms in $D_{\xi }^{p}\left\langle f(x),\eta \right\rangle ,$ 
$1\leq p\leq k,$ except possibly 
\begin{equation}
\left( \ln \left\vert \ln \left\vert x\right\vert \right\vert \right)
\left\langle D_{\xi }^{k}\varphi _{\cdot }(x),\eta \right\rangle
\label{k-th term}
\end{equation}%
tend to $0$ as $x\rightarrow 0.$ We only need to prove that the function in (%
\ref{k-th term}) also tends to $0$ as $x\rightarrow 0.$ By expanding $\left(
x_{1}+t\xi _{1}\right) ^{\alpha _{1}}\left( x_{2}+t\xi _{2}\right) ^{\alpha
_{2}}...\left( x_{n}+t\xi _{n}\right) ^{\alpha _{n}}$\ binomially,\ we can
write%
\begin{equation*}
\varphi _{\cdot }(x+t\xi ):=\varphi _{\cdot }(x)+P(x,\xi ,t)+\varphi _{\cdot
}(\xi )t^{k}.
\end{equation*}%
But since $(\xi ,\eta )\in \Lambda ,$%
\begin{equation*}
\left\langle D_{\xi }^{k}\varphi _{\cdot }(x),\eta \right\rangle
=k!\left\langle \varphi _{\cdot }(\xi ),\eta \right\rangle =0.
\end{equation*}

It follows that $\left\vert D_{\xi }^{p}\left\langle f(0),\eta \right\rangle
\right\vert =0$\ for $p\leq k.$\ Thus, $f\in C_{\xi }^{k}$\ for all $(\xi
,\eta )\in \Lambda ,$\ but $f\notin C^{k}$.
\end{proof}

\begin{remark}
(cf. \cite{korevaar}) Suppose (i) is satisfied for all $k\geq 0.$ It would
be of interest to know whether there exists a constant $\rho =\rho (\Lambda
),$ depending only on some approperiate notion of capacity of $\Lambda ,$ so
that (\ref{B}) is satisfied with $B=\left( \rho (\Lambda )\right) ^{-k}$ for
all $f$ and and all $k.$
\end{remark}

\begin{remark}
Suppose $\Lambda $\ satisfies (i) or (ii). The proof of Theorem \ref{main}
shows that if $f$\ is continuous and $D_{\xi }^{k}\left\langle f,\eta
\right\rangle =0$,$\forall (\xi ,\eta )\in \Lambda ,$\ then $f$\ is a
polynomial. The assumption of continuity of $f$\ is not necessary but our
proof is valid only if $f$\ is continuous. See \cite{Boman-partial}.
\end{remark}

\begin{remark}
If $\Lambda $\ satisfies (i), then $\Lambda $\ contains at least $mk_{n}$\
elements. Furthermore, if (i) holds for $k$\ then (i) also holds for all $%
j\leq k.$ Suppose there exists $\Phi \in B_{j},j<k$\ such that $\left. \Phi
\right\vert _{\Lambda }\equiv 0$\ but $\Phi \not\equiv 0.$\ Then, $\xi
_{1}^{k-j}\Phi \in \mathcal{B}_{k}$, $\left. \xi _{1}^{k-j}\Phi \right\vert
_{\Lambda }\equiv 0$\ but this is a contradiction.
\end{remark}

\ Let $\left\{ M_{k}\right\} _{k=0}^{\infty },$ be a sequence of nonnegative
numbers. For $h>0$\ and $K\subset \subset 
%TCIMACRO{\U{211d} }%
%BeginExpansion
\mathbb{R}
%EndExpansion
^{n}$ define the seminorm on $C^{\infty }\left( 
%TCIMACRO{\U{211d} }%
%BeginExpansion
\mathbb{R}
%EndExpansion
^{n}\right) ,$%
\begin{equation*}
p_{h,K}(f)=\sup_{\alpha \in 
%TCIMACRO{\U{2124} }%
%BeginExpansion
\mathbb{Z}
%EndExpansion
_{+}^{n}}\sup_{x\in K}\frac{\left\vert \partial ^{\alpha }f(x)\right\vert }{%
h^{|\alpha |}M_{|\alpha |}}.
\end{equation*}%
The spaces

\begin{equation*}
C\left\{ M_{k}\right\} =\left\{ f\in C^{\infty }(%
%TCIMACRO{\U{211d} }%
%BeginExpansion
\mathbb{R}
%EndExpansion
^{n}):\forall K\subset \subset 
%TCIMACRO{\U{211d} }%
%BeginExpansion
\mathbb{R}
%EndExpansion
^{n},\exists h>0,\text{ s.t. }p_{h,K}(f)<\infty \right\}
\end{equation*}%
and%
\begin{equation*}
\text{ }C\left( M_{k}\right) =\left\{ f\in C^{\infty }(%
%TCIMACRO{\U{211d} }%
%BeginExpansion
\mathbb{R}
%EndExpansion
^{n}):p_{h,K}(f)<\infty ,\forall K\subset \subset 
%TCIMACRO{\U{211d} }%
%BeginExpansion
\mathbb{R}
%EndExpansion
^{n},\forall h>0\right\}
\end{equation*}%
\ are called the Carleman and Beurling classes, respectively. The classes $%
C\left\{ (k!)^{\nu }\right\} ,$ $\nu >1$,\ known as Gevrey classes, are
especially important in partial differential equations and harmonic
analysis. The class $C\left\{ k!\right\} $ is precisely the class of real
analytic functions.

We assume that 
\begin{equation}
M_{0}=1\text{ and }M_{k}\geq k!,\forall k;  \label{M_k>k!}
\end{equation}%
\begin{equation}
M_{k}^{1/k}\text{ is strictly increasing;}  \label{M_k  inc}
\end{equation}%
\qquad 
\begin{equation}
\exists C>0\text{ such that }M_{k+1}\leq C^{k}M_{k},\text{ }\forall k.
\label{M_k diff}
\end{equation}%
These conditions insure that the classes $C\{M_{k}\}$ and $C(M_{k})$ are
nontrivial and are closed under product and differentiation of functions.%
\textrm{\ }For more properties of these spaces, see \cite{Hormander}, \cite%
{n4} and references there.

It is well known that $f\in C^{\infty }(%
%TCIMACRO{\U{211d} }%
%BeginExpansion
\mathbb{R}
%EndExpansion
^{n})$\ if and only if $\sup_{\xi \in 
%TCIMACRO{\U{211d} }%
%BeginExpansion
\mathbb{R}
%EndExpansion
^{n}}\left\vert \xi \right\vert ^{j}|\widehat{\chi f}(\xi )|<\infty ,\forall
\chi \in C_{c}^{\infty }(%
%TCIMACRO{\U{211d} }%
%BeginExpansion
\mathbb{R}
%EndExpansion
^{n}),j\geq 1.$ A similar characterization is also available for $C\{M_{k}\}$%
\ (see \cite{Hormander}) a routine modification of which yields an analogous
characterization of $C(M_{k})$.

Let $r>0$. Choose a sequence of cut-off functions $\chi _{(j)}\in
C_{c}^{\infty },j=1,2,...,$\ such that $\chi _{(j)}(x)=1$\ if $|x-x_{0}|<r$, 
$\chi _{(j)}(x)=0$\ if $|x-x_{0}|>3r$ and 
\begin{equation*}
\left\vert \partial ^{\alpha }\chi _{(j)}(x)\right\vert \leq \left(
C_{1}j\right) ^{\left\vert \alpha \right\vert },\forall j,\forall \left\vert
\alpha \right\vert \leq j,\forall x,
\end{equation*}%
where the constant $C_{1}$\ is independent of $j.$

Then $f\in C\left\{ M_{k}\right\} $\ (\emph{resp.} $C\left( M_{k}\right) $)
in a neighborhood of $x_{0}\in 
%TCIMACRO{\U{211d} }%
%BeginExpansion
\mathbb{R}
%EndExpansion
^{n}$\ if and only if there exists a constant $\hbar >0$\ (\emph{resp.} for
every $\hbar >0$) 
\begin{equation*}
\sup_{\xi \in 
%TCIMACRO{\U{211d} }%
%BeginExpansion
\mathbb{R}
%EndExpansion
^{n}}\sup_{j\geq 1}\hbar ^{-j}M_{j}^{-1}|\xi |^{j}|\widehat{f\chi _{(j)}}%
(\xi )|<\infty .
\end{equation*}

Call a subset $\Lambda \subset 
%TCIMACRO{\U{211d} }%
%BeginExpansion
\mathbb{R}
%EndExpansion
^{n}\times 
%TCIMACRO{\U{211d} }%
%BeginExpansion
\mathbb{R}
%EndExpansion
^{m}$\ a determining set for bilinear forms of rank 1 if there is no nonzero
bilinear form $\varphi (\xi ,\eta ),\xi \in 
%TCIMACRO{\U{211d} }%
%BeginExpansion
\mathbb{R}
%EndExpansion
^{n},\eta \in 
%TCIMACRO{\U{211d} }%
%BeginExpansion
\mathbb{R}
%EndExpansion
^{m}$\ of rank 1 such that $\varphi (\xi ,\eta )=0$\ for all $(\xi ,\eta
)\in \Lambda .$

Clearly $\Lambda $\ is a determining set for bilinear forms of rank 1 if and
only if 
\begin{equation*}
\left\langle u,\xi \right\rangle \left\langle v,\eta \right\rangle
=0,\forall (\xi ,\eta )\in \Lambda \Rightarrow |u||v|=0
\end{equation*}%
(here $\left\langle u,\xi \right\rangle $\ and $\left\langle v,\eta
\right\rangle $\ are dot products on $%
%TCIMACRO{\U{211d} }%
%BeginExpansion
\mathbb{R}
%EndExpansion
^{n}$\ and $%
%TCIMACRO{\U{211d} }%
%BeginExpansion
\mathbb{R}
%EndExpansion
^{m}$, respectively), or equivalently,\ 
\begin{equation*}
\cap _{(\xi ,\eta )\in \Lambda }\{(u,v)\in 
%TCIMACRO{\U{211d} }%
%BeginExpansion
\mathbb{R}
%EndExpansion
^{n}\times 
%TCIMACRO{\U{211d} }%
%BeginExpansion
\mathbb{R}
%EndExpansion
^{m}:<u,\xi ><v,\eta >=0\}=\left( 
%TCIMACRO{\U{211d} }%
%BeginExpansion
\mathbb{R}
%EndExpansion
^{n}\times 0\right) \cup (0\times 
%TCIMACRO{\U{211d} }%
%BeginExpansion
\mathbb{R}
%EndExpansion
^{m}).
\end{equation*}%
Since $%
%TCIMACRO{\U{211d} }%
%BeginExpansion
\mathbb{R}
%EndExpansion
\lbrack u,v]$\ is a Noetherian ring, $\Lambda $\ contains a finite subset $%
\Lambda ^{\prime }$\ such that the sets $\{<u,\xi ><v,\eta >:(\xi ,\eta )\in
\Lambda \}$\ and $\{<u,\xi ><v,\eta >:(\xi ,\eta )\in \Lambda ^{\prime }\}$\
generate the same ideal in $%
%TCIMACRO{\U{211d} }%
%BeginExpansion
\mathbb{R}
%EndExpansion
\lbrack u,v]$\ and thus define the same varieties:%
\begin{eqnarray*}
\cap _{(\xi ,\eta )\in \Lambda }\{(u,v) &\in &%
%TCIMACRO{\U{211d} }%
%BeginExpansion
\mathbb{R}
%EndExpansion
^{n}\times 
%TCIMACRO{\U{211d} }%
%BeginExpansion
\mathbb{R}
%EndExpansion
^{m}:<u,\xi ><v,\eta >=0\} \\
&=&\cap _{(\xi ,\eta )\in \Lambda ^{\prime }}\{(u,v)\in 
%TCIMACRO{\U{211d} }%
%BeginExpansion
\mathbb{R}
%EndExpansion
^{n}\times 
%TCIMACRO{\U{211d} }%
%BeginExpansion
\mathbb{R}
%EndExpansion
^{m}:<u,\xi ><v,\eta >=0\}.
\end{eqnarray*}%
Thus, any determining set for bilinear forms of rank 1 contains a finite
determining set for bilinear forms of rank 1.

Let $C\left\{ M_{k}\right\} \left( \xi \right) $\ (\emph{resp.} $C\left(
M_{k}\right) (\xi )$) denote the set of all $f\in C_{\xi }^{\infty }(%
%TCIMACRO{\U{211d} }%
%BeginExpansion
\mathbb{R}
%EndExpansion
^{n})$\ such that for every subset $K\subset \subset 
%TCIMACRO{\U{211d} }%
%BeginExpansion
\mathbb{R}
%EndExpansion
^{n},$ $\sup_{j,x\in K}\left\vert D_{\xi }^{j}f(x)\right\vert \hbar
^{-j}M_{j}^{-1}<\infty ,\forall j,$ for some $\hbar >0$\ (\emph{resp.}\ for
every $\hbar >0$).

\begin{theorem}
Let $\{M_{k}\}_{k=0}^{\infty }$ be a sequence of nonnegative numbers
satisfying the conditions ( \ref{M_k>k!}), (\ref{M_k inc}) and (\ref{M_k
diff}). The following statements are equivalent.

(i) $\Lambda $\ is a determining set for bilinear forms of rank 1;

(ii) for any locally bounded map $f:%
%TCIMACRO{\U{211d} }%
%BeginExpansion
\mathbb{R}
%EndExpansion
^{n}\rightarrow 
%TCIMACRO{\U{211d} }%
%BeginExpansion
\mathbb{R}
%EndExpansion
^{m}$,\ 
\begin{equation*}
\left\langle \eta ,f\right\rangle \in C\left\{ M_{k}\right\} \left( \xi
\right) ,\forall (\eta ,\xi )\in \Lambda \Rightarrow f\in C\left\{
M_{k}\right\} ;
\end{equation*}%
(iii) for any locally bounded map $f:%
%TCIMACRO{\U{211d} }%
%BeginExpansion
\mathbb{R}
%EndExpansion
^{n}\rightarrow 
%TCIMACRO{\U{211d} }%
%BeginExpansion
\mathbb{R}
%EndExpansion
^{m}$, 
\begin{equation*}
\left\langle \eta ,f\right\rangle \in C\left( M_{k}\right) \left( \xi
\right) ,\forall (\eta ,\xi )\in \Lambda \Rightarrow f\in C\left(
M_{k}\right) ;
\end{equation*}%
(iv) for any locally bounded map $f:%
%TCIMACRO{\U{211d} }%
%BeginExpansion
\mathbb{R}
%EndExpansion
^{n}\rightarrow 
%TCIMACRO{\U{211d} }%
%BeginExpansion
\mathbb{R}
%EndExpansion
^{m}$,%
\begin{equation*}
\left\langle \eta ,f\right\rangle \in C^{\infty }\left( \xi \right) ,\forall
(\eta ,\xi )\in \Lambda \Rightarrow f\in C^{\infty }.
\end{equation*}
\end{theorem}

\begin{proof}
(cf. Theorem 4 in \cite{Boman-partial}) Assume (i) holds. By the remark
above, by replacing $\Lambda $ by a subset, if necessary, we may assume $%
\Lambda $ is finite. Suppose for every $(\eta ,\xi )\in \Lambda ,$\ $%
\left\langle \eta ,f\right\rangle \in C\left\{ M_{k}\right\} \left( \xi
\right) $\ (resp. $\left\langle \eta ,f\right\rangle \in C\left(
M_{k}\right) \left( \xi \right) $). Now for a suitable function $f,$ 
\begin{equation*}
\left\langle \xi ,z\right\rangle \widehat{\left\langle \eta ,f\right\rangle }%
(z)=\left\langle \xi ,z\right\rangle \left\langle \eta ,\widehat{f}%
(z)\right\rangle =\left\langle \eta ,i\int \left[ \left\langle \xi ,\partial
_{x}\right\rangle e^{-i\left\langle x,z\right\rangle }\right]
f(x)dx\right\rangle
\end{equation*}%
\begin{equation*}
=\left\langle \eta ,-i\int e^{-i\left\langle x,z\right\rangle }\left\langle
\xi ,\partial _{x}f\right\rangle (x)dx\right\rangle =\left\langle \eta
,-i\int e^{-i\left\langle x,z\right\rangle }D_{\xi }f(x)dx\right\rangle .
\end{equation*}%
Let $g_{(j)}:=f\chi _{(j)}\in C\left\{ M_{k}\right\} $\ near a fixed point $%
x_{0}.$\ Assume, without loss of generality, $x_{0}=0$. By assumption, for
all $(\xi ,\eta )\in \Lambda $ there exist constants $C=C_{\xi \eta }$ and $%
\hbar =\hbar _{\xi \eta }>0$\ (\emph{resp.} for all $(\xi ,\eta )\in \Lambda 
$ and for all $\hbar >0$ there exists a constant $C=C_{\xi \eta ,\hbar }$ )
such that \ 
\begin{equation*}
\left\vert \widehat{\left\langle \eta ,g_{(j)}\right\rangle }(\zeta
)\right\vert \left\vert \left\langle \xi ,\zeta \right\rangle \right\vert
^{j}=\left\vert \left\langle \eta ,\widehat{g_{(j)}}(\zeta )\right\rangle
\right\vert \left\vert \left\langle \xi ,\zeta \right\rangle \right\vert
^{j}\leq C\hbar ^{j}M_{j},\forall \left( \xi ,\eta \right) \in \Lambda
,\zeta \in 
%TCIMACRO{\U{211d} }%
%BeginExpansion
\mathbb{R}
%EndExpansion
^{n},j\in 
%TCIMACRO{\U{2124} }%
%BeginExpansion
\mathbb{Z}
%EndExpansion
_{+}.
\end{equation*}

The function 
\begin{equation}
%TCIMACRO{\U{211d} }%
%BeginExpansion
\mathbb{R}
%EndExpansion
^{n}\times 
%TCIMACRO{\U{211d} }%
%BeginExpansion
\mathbb{R}
%EndExpansion
^{m}\ni (u,v)\rightarrow \sum_{(\xi ,\eta )\in \Lambda }\left\vert
\left\langle \eta ,v\right\rangle \right\vert \left\vert \left\langle \xi
,u\right\rangle \right\vert ^{l},  \label{lowerbound}
\end{equation}%
is homogeneous of degree $1$\ in $v$, of homogeneous degree $l$\ in $u.$\
Since none of the terms $\left\vert \left\langle \eta ,v\right\rangle
\right\vert \left\vert \left\langle \xi ,u\right\rangle \right\vert $\ can
vanish on all of $\Lambda ,$\ the function in (\ref{lowerbound}) has a
positive minimum on the compact set $\{(u,v):|u|=1,|v|=1\}.$\ Thus, there is
an $\varepsilon >0$ such that 
\begin{equation*}
\sum_{(\xi ,\eta )\in \Lambda }\left\vert \left\langle \eta ,v\right\rangle
\right\vert \left\vert \left\langle \xi ,u\right\rangle \right\vert ^{l}\geq
\varepsilon |v||u|^{l},
\end{equation*}%
(see Lemma 1 \cite{Boman-partial}). Applying this to $u=\zeta ,v=\widehat{%
g_{(j)}}(\zeta ),$\ we get 
\begin{equation*}
\left\vert \widehat{g_{(j)}}(\zeta )\right\vert \left\vert \zeta \right\vert
^{l}\leq \varepsilon ^{-1}\sum_{(\xi ,\eta )\in \Lambda }\left\vert
\left\langle \eta ,\widehat{g_{(j)}}(\zeta )\right\rangle \right\vert
\left\vert \left\langle \xi ,\zeta \right\rangle \right\vert ^{l}\leq C\hbar
^{j}M_{j},
\end{equation*}%
where $\hbar =\max_{(\xi ,\eta )\in \Lambda }\hbar _{\xi \eta }$\ (\emph{resp%
}. for all $\hbar >0$) and $C=\varepsilon ^{-1}\sum_{(\xi ,\eta )\in \Lambda
}C_{\xi \eta }.$\ Thus (ii) and (iii) hold. By setting $\hbar =1$\ and $%
M_{j}=1,\forall j,$\ in the above argument, it is clear that (iii) holds as
well.

Conversely if $\Lambda $\ is not a determinant set for bilinear forms of
rank 1, there exist $u\neq 0$\ and $v\neq 0$\ such that 
\begin{equation*}
\left\langle u,\xi \right\rangle \left\langle v,\eta \right\rangle
=0,\forall \left( \xi ,\eta \right) \in \Lambda .
\end{equation*}%
Let $h:%
%TCIMACRO{\U{211d} }%
%BeginExpansion
\mathbb{R}
%EndExpansion
\rightarrow 
%TCIMACRO{\U{211d} }%
%BeginExpansion
\mathbb{R}
%EndExpansion
$\ be an arbitrary continuous function. Let $f:%
%TCIMACRO{\U{211d} }%
%BeginExpansion
\mathbb{R}
%EndExpansion
^{n}\rightarrow 
%TCIMACRO{\U{211d} }%
%BeginExpansion
\mathbb{R}
%EndExpansion
^{m}$\ be defined as $f(z)=h\left( \left\langle u,z\right\rangle \right)
\cdot v.$\ Then 
\begin{equation*}
\left. \left( \frac{d}{dt}\left\langle \eta ,f(z+t\xi )\right\rangle \right)
\right\vert _{t=0}=\left. \left\langle \eta ,v\right\rangle \left\langle
u,\xi \right\rangle h^{\prime }\left( \left\langle u,z+t\xi \right\rangle
\right) \right\vert _{t=0}\equiv 0
\end{equation*}%
Thus $\left\langle \eta ,f\right\rangle \in C(M_{k})\left( \xi \right)
\subset C\{M_{k}\}\left( \xi \right) \subset C^{\infty }\left( \xi \right)
,\forall (\xi ,\eta )\in \Lambda $\ but $f$\ need not be even differentiable.
\end{proof}

\end{document}